          \documentclass{cmslatex}
\usepackage[paperwidth=7in, paperheight=10in, margin=.875in]{geometry}
\usepackage{amsmath}
\usepackage{setspace,graphicx}     
\usepackage{amsfonts}
\usepackage{latexsym}
\usepackage{upgreek}
\usepackage{amssymb}

\usepackage{amsthm}
\usepackage{setspace}
\usepackage{fancyhdr}
\usepackage{mathrsfs}
\usepackage{framed}
\usepackage{bm}
\usepackage{perpage} 
\usepackage{bbm}
\usepackage{dsfont}
\usepackage{comment}
\usepackage{color}
\definecolor{darkgreen}{rgb}{0,0.7,0}

\definecolor{darkblue}{rgb}{0.2,0,0.7}

\newcommand{\R}{\mathbb{R}}
\newcommand{\rp}{\mathbb{P}}

\newcommand{\N}{\mathbb{N}}

\newcommand{\dkl}{D_{\mbox {\tiny{\rm KL}}}}
\newcommand{\dtv}{d_{\mbox {\tiny{\rm TV}}}}

\newcommand{\E}{\mathbb{E}}

\newcommand{\envelope}{(\raisebox{-.5pt}{\scalebox{1.45}{\Letter}}\kern-1.7pt)}

\renewcommand{\theequation}{\arabic{section}.\arabic{equation}}

            \newtheorem{thm}{Theorem}[section]
          
          \newtheorem{lem}[thm]{Lemma}
          
          \newtheorem{remark}[thm]{Remark}
		\newtheorem{assumption}[thm]{Assumption}

          \newtheorem{example}[thm]{Example}

\begin{document}
          \title{Gaussian Approximations of Small Noise Diffusions in Kullback-Leibler Divergence\thanks{}DS-A is supported by EPSRC as part of the MASDOC DTC at the University of Warwick with grant No. EP/H023364/1. The work of AMS is supported by DARPA, EPSRC and ONR.}


          \author{Daniel Sanz-Alonso\thanks{Mathematics Institute, University of
Warwick, Coventry CV4 7AL, United Kingdom. d.sanz-alonso@warwick.ac.uk}
         \and{Andrew M. Stuart\thanks{ Mathematics Institute, University of
Warwick Coventry CV4 7AL, United Kingdom, a.m.stuart@warwick.ac.uk}}}





         \pagestyle{myheadings} \markboth{Gaussian Approximations of Small Noise Diffusions}{D. Sanz-Alonso and A. M. Stuart} \maketitle

          \begin{abstract}
We study Gaussian approximations to the distribution of a diffusion. 
The approximations are easy to compute: they are defined by two simple ordinary 
differential equations for the mean and the covariance. Time correlations can also be computed via solution of
a linear stochastic differential equation. We show, using the Kullback-Leibler divergence, that the approximations are accurate in the small noise regime. 
An analogous discrete time setting is also studied.
The results provide both theoretical support for the use of Gaussian processes in the approximation of diffusions, and methodological guidance in the construction of Gaussian approximations in applications.
          \end{abstract}
\begin{keywords}
Gaussian approximations, diffusion processes, small noise, Kullback-Leibler divergence
\end{keywords}

 \begin{AMS} 28C20, 60H10, 65L05, 60G15
\end{AMS}
          \section{Introduction}\label{intro}
Consider the stochastic differential equation (SDE) 
\begin{equation}\label{sde}
dv^\epsilon (t)= f\bigl(v^\epsilon(t)\bigr) dt + \sqrt{\epsilon\Sigma}\, dW(t), \quad t\in [0,T], \quad v^\epsilon(0)\sim \mu_0^\epsilon.
\end{equation}
Here, and throughout, the time-horizon $T>0$ is finite and fixed, the drift $f:\, \R^D \to \R^D$ is nonlinear in the case of interest, $\Sigma$ is a positive-definite matrix, and $dW$ is a standard $D$-dimensional Wiener process. Assume that the paths $v^\epsilon$ are continuous, and let $\mu^\epsilon$ be the law of $v^\epsilon$ in $C([0,T],\R^D)$. We show that in the small $\epsilon$ regime ---corresponding to small diffusion coefficient and little uncertainty in the initial condition---  $\mu^\epsilon$ can be accurately approximated by a Gaussian measure $\nu^\epsilon$ in $C([0,T],\R^D)$, whose mean function and marginal covariances satisfy simple ordinary differential equations (ODEs). Precisely, we show that the Kullback-Leibler divergence $\dkl(\nu^\epsilon\|\mu^\epsilon)$ is of order $\epsilon.$  We investigate the effect that numerically approximating these ODEs has on the approximation of $\mu^\epsilon.$ Finally, we construct ---and show the accuracy of--- Gaussian approximations to Euler-Maruyama discretizations of \eqref{sde}.

We aim to provide a rigorous justification of Gaussian process approximations of small noise diffusions. The use of Gaussian processes is now pervasive in applications \cite{williams}. However, Gaussian assumptions are often used for algorithmic and mathematical convenience without much theoretical support. An example of this is in the field of data assimilation, where many algorithms invoke Gaussian approximations in order to apply Kalman formulae \cite{dataassimilationbook}. The paper \cite{archambeau} suggests a variational approach to computing Gaussian process approximations to posterior measures arising from discretely observed diffusions. The algorithm proposed aims to minimize the Kullback-Leibler divergence between the Gaussian process and the posterior of interest. In \cite{weber2}, the authors study the well-posedness of the abstract problem of finding the best (in the Kullback-Leibler sense) Gaussian approximation $\nu$ to a measure $\mu$. They use a nonconstructive calculus of variation approach. In the complementary paper \cite{weber1}, the authors propose an algorithm that allows 
numerical computation of such best approximation. The paper \cite{yulong} seeks the best Kullback-Leibler Gaussian approximation to the distribution of conditioned diffusions arising in molecular dynamics. Here our interest is not in finding the best Gaussian approximation, but in constructing Gaussians that can be proved to approximate accurately the distribution of a small noise diffusion. The Gaussian approximations are built using certain linearizations of the drift. Crucially, the time-marginals of our approximations can be found by solving (perhaps numerically) two simple ODEs, which can be advantageous over numerically solving the full Fokker-Planck equation. Our computable Gaussian approximations and error bounds are of potential interest in filtering signals arising from small noise diffusions, for instance within the methodology  in \cite{lawfokker}.

The paper is organized as follows. Subsection \ref{ssec:odesintro} sets the mathematical framework, and introduces the ODEs that will be used to define the means and covariances of the time-marginals of our approximations.  Subsection \ref{ssec:kullbackleiblerbackground} gives some background on the Kullback-Leibler divergence. Subsection \ref{ssec:largedeviationsjustification} motivates the choice of the ODEs. Section \ref{sec:cts-cts} studies the approximation of the SDE using a Gaussian whose time-marginals  mean and covariance satisfy the ODEs 
given in Subsection \ref{ssec:odesintro}. In Section \ref{sec:cts-disc} we investigate the effect that numerically solving the ODEs has on the quality of the resulting Gaussian approximation. Finally, Section \ref{sec:disc-disc} studies the approximation of Euler discretizations of \eqref{sde} via Gaussians with means and covariances defined through discretization of the ODEs.

\paragraph{Notation} We denote by $|\cdot |$ the standard inner product in $\R^D.$ For positive-definite $P\in \R^{D\times D}$ we denote $| \cdot|_P:= |P^{-1/2} \cdot|.$  Subscripts in path measures will denote time-marginals. Thus, for a measure $\mu$ in $C([0,T],\R^D)$ and $0\le t\le T,$  we denote by $\mu_t$ the push-forward measure in $\R^D$ via the evaluation map $E_t(h) := h(t),$  $h\in C([0,T],\R^D)$. In the discrete setting, given a path measure $\mu_{0:k}$ in $\R^{D(k+1)}$ and $0\le j \le k,$  $\mu_j$ will denote the push-forward measure in $\R^D$ via the evaluation map $E_j:(v_d)_{d=1}^{D(k+1)}\in \R^{D(k+1)} \mapsto (v_{j+d})_{d=1}^D \in \R^D$. A subscript $d$ in a function or vector denotes, as usual, component and coordinate. Superscripts highlight parameters of interest in the analysis.

\subsection{Mathematical Framework. ODEs for Mean and Covariance}\label{ssec:odesintro}
We work under the following assumption 
(recalling that $T$ is fixed and finite):
\begin{assumption}\label{assumptionode}
$f\in C^2(\R^D,\R^D)$ and  there exist $c>0$ and $s\in \N$ such that, for $1\le d\le D,$ and every multi-index $\alpha$ with $|\alpha|=2,$  $|\partial^{\alpha} f_{d}(u)|\le c(1+|u|^s).$ Moreover, for any $m_0 \in \R^D,$  the ODE
\begin{equation}\label{ODE}
\frac{dm}{dt} = f(m), \quad t\in [0,T],  \quad m(0) = m_0\in \R^D,
\end{equation}
has a unique solution.
\end{assumption}

Equation \eqref{ODE}, or discretizations thereof, will be used to define the mean of our Gaussian approximations. The  covariance of the time-marginals will be given by $C^\epsilon=\epsilon C$ where
\begin{equation}\label{ODEcov}
\frac{dC}{dt} = Df(m)C + C Df(m)^T + \Sigma, \quad t\in [0,T],  \quad C(0) = 
C_0\in \R^{D\times D},
\end{equation}
or discretizations thereof, and where $C_0$ is a given positive semi-definite matrix ---see Remark \ref{initialcovariance}.
This equation has a unique solution on $[0,T]$ under Assumption \ref{assumptionode}.

It is intuitively clear that, for small $\epsilon,$ the solution $v^\epsilon=\{v^\epsilon (t)\}_{0\le t \le T}$ of \eqref{sde} is close to the deterministic solution $m=\{m(t)\}_{0\le t \le T}$ of \eqref{ODE}. In Section \ref{sec:cts-cts} we show that, for small $\epsilon,$ the law of $v^\epsilon(t)$ is approximately Gaussian with mean $m(t)$ and covariance $C^\epsilon(t)=\epsilon C(t)$. More precisely, the error in this approximation, measured in the Kullback-Leibler divergence, is of order $\epsilon.$ Indeed we show more: the law  $\mu^{\epsilon}$ of the paths $v^\epsilon$ can be approximated by a Gaussian measure $\nu^{\epsilon}$ in $C([0,T],\R^D).$ Subsection \ref{ssec:largedeviationsjustification} shows that the choice of $m$ as the mean of the Gaussian approximation is crucial. If any other function is chosen as mean, the Kullback-Leibler error remains 
of order $1$.

\begin{remark}\label{initialcovariance}
Our results will cover two possible initializations of the SDE \eqref{sde}. First, deterministic initial condition $v^\epsilon(0) = v_0$, corresponding to $\mu_0^\epsilon = \delta_{v_{0}}$. In this case \eqref{ODE} and \eqref{ODEcov} should be initialized with $m_0=v_0$ and $C_0=0,$  corresponding to the same Dirac measure. Second, where $\mu_0^{\epsilon}$ has positive Lebesgue density in $\R^D$, and $\dkl\Bigr(N(m_0, \epsilon C_0)\|\mu_0^\epsilon\Bigl)$ is of order $\epsilon.$
\end{remark}

\begin{remark}\label{rem:disc-cov-bound}
In Sections \ref{sec:cts-disc} and \ref{sec:disc-disc} we discretize \eqref{ODE} and \eqref{ODEcov} with step-size $\Delta t.$  Assumption \ref{assumptionode} and the convergence of the methods employed imply that, for all sufficiently small $\Delta t,$ the discretized means and covariances are uniformly bounded. 
That is, there is $M>0$ such that, for all sufficiently small $\Delta t$, 
$|m_k^{\Delta t}|, \,|C_k^{\Delta t}|\le M$, provided that $0 \le k\Delta t \le T.$ 
\end{remark}

		\subsection{Background on Kullback-Leibler Divergence}\label{ssec:kullbackleiblerbackground}
Let $\nu$ and $\mu$ be two probability measures on a measurable space $(X, {\cal X})$. The Kullback-Leibler divergence (also known as relative entropy) of $\nu$ with respect to $\mu$ is given by 
\begin{equation*}
\dkl(\nu\|\mu):= \E^{\nu} \log\Bigl(\frac{d\nu}{d\mu}\Bigr)
\end{equation*}
if $\nu$ is absolutely continuous with respect to $\mu$, denoted $\nu\ll \mu,$ and $\dkl(\nu\|\mu)=\infty$ otherwise. 
The Kullback-Leibler divergence satisfies  $\dkl(\nu\|\mu)\ge 0,$ but it is not a metric on the space of probability measures since it may not be finite, it is not symmetric, and it does not satisfy the triangle inequality. However, it does quantify the proximity of the measures $\nu$ and $\mu.$ For instance, it provides an upper bound on the total variation distance
\begin{equation*}
\dtv(\nu,\mu):= \sup {\big \{ }|\nu(A) - \mu(A)| : A\in {\cal X} \big{ \} }
\end{equation*}
 via Pinsker's inequality \cite{pinsker1960information}
\begin{equation}\label{pinskerine}
\dtv(\nu,\mu)\le \dkl(\nu\|\mu)^{1/2}.
\end{equation}%

We will use the chain rule for Kullback-Leibler divergence \cite{dupuis2011weak}, which is well known in information theory. We recall it in the next lemma.
The result underlies the proof of Lemma \ref{lemmageneralbounddiscretetime},
which is analogous to the continuous time result
Lemma \ref{lemmaconttime}.

\begin{lem}\label{lemmakl}
Let $X$ and $Y$ be Polish spaces and $\nu$ and $\mu$ be probability measures on the measurable space $(X\times Y, {\cal  F}).$ Denote by $\nu_x$ and $\mu_x$ the first marginals of $\nu$ and $\mu$ and let $\nu(dy|x)$ and $\mu(dy|x)$ be stochastic kernels on $Y$ given $X$ for which, for $A\times B \in {\cal F},$ 
\begin{equation*}
\nu(A\times B) = \int_A \nu(B|x) \nu_x(dx),  \quad \mu(A\times B) = \int_A \mu(B|x) \mu_x(dx).
\end{equation*}
Then,
\begin{equation}\label{conditional}
\dkl(\nu\| \mu) = \dkl(\nu_x\|  \mu_x) + \E^{\nu_x} \dkl\bigl(\nu(\cdot|x) \|  \mu(\cdot|x)\bigr).
\end{equation}
\end{lem}

The non-negativity of the Kullback-Leibler divergence and \eqref{conditional} imply that any two corresponding marginals are closer in Kullback-Leibler than the full measures, i.e.
\begin{equation}\label{marginalsbound}
\dkl(\nu_x\|  \mu_x) \le \dkl(\nu\|\mu).
\end{equation}

The chain rule is  powerful when studying the approximation, by a Gaussian measure $\nu$, of a non-Gaussian measure $\mu$ for which $\mu(\cdot |x)$ is Gaussian: the structure in \eqref{conditional} allows exploitation of the Gaussianity of the kernels $\nu(\cdot|x)$ and $\mu(\cdot|x).$ This will become apparent in Section \ref{sec:disc-disc}, Lemma \ref{lemmageneralbounddiscretetime}.

\subsection{Large Deviations and the Choice of Mean}\label{ssec:largedeviationsjustification}
Small noise diffusions have been extensively studied using large deviations. An early and fundamental result in the theory can be found in the first edition of \cite{freidlinwentzell}, where it was shown that the collection $\{v^\epsilon, \epsilon\in (0,1)\}$ defined by \eqref{sde} with initial condition $v^\epsilon(0)=  v_0$ satisfies a large deviation principle on $C([0,T], \R^D)$ with rate function 

\begin{equation}\label{ratefunction}
I(\varphi) := \inf_{u\in U_\varphi} \bigg{ \{ } \frac12 \int_0^T |u(t)|^2 \, dt \bigg{\}},
\end{equation}
where 
\begin{equation*}
U_{\varphi} := {\biggl \{}  u \in L^2([0,T],\R^D) : \varphi(t) = v_0 + \int_0^T f\bigl(\varphi(s)\bigr) \,ds + \sqrt{ \Sigma }\int_0^t u(s) \,ds  {\biggr \}},
\end{equation*}
for absolutely continuous $\varphi$ with $\varphi(0)=v_0$, $U_{\varphi}=\emptyset$
 for all other $\varphi \in C([0,T], \R^D)$, and the infimum over the empty set in \eqref{ratefunction} is taken to be $\infty.$ 
 
It follows from the definition of the rate function $I$ in \eqref{ratefunction} that $I(\varphi) = 0$ iff the zero function $u\equiv 0$ is in  $U_{\varphi}$. This holds iff $\varphi$ is the solution to the ODE \eqref{ODE} with initial condition $\varphi(0) = v_0.$ 
In a very rough sense, this implies that the probability of $v^\epsilon$ lying in a small tube centered around any function other than the solution to the ODE decays at least exponentially as $\epsilon \to 0.$  That is, for small $\delta, \epsilon>0,$ we have, as $\epsilon \to 0,$
\begin{equation*}
\rp \bigl( v^\epsilon \in B_\delta (\varphi)\bigr) \approx \exp \biggl(-\frac{I(\varphi)}{\epsilon}\biggr)
\end{equation*}
where $B_\delta(\varphi)$ denotes the ball (in the supremum norm) of radius $\delta>0$  and center $\varphi\in C([0,T], \R^D).$ 

Combining the large deviation result and the Markov inequality we obtain the following:
\begin{lem}\label{lemmaldp}
Let $\mu^\epsilon$ be the law in $C([0,T],\R^D)$ of $v^\epsilon$  given by \eqref{sde} with initial condition $v^{\epsilon}(0)= v_0$.  Consider a Gaussian measure $\nu^\epsilon$ in $C([0,T],\R^D)$ with mean $m\in C([0,T],\R^D)$ and time-marginal covariances $C^\epsilon(t) = \epsilon C(t)$. Then, unless $m$ is the solution to \eqref{ODE} with initial condition $m(0)=v_0$ we have, as $\epsilon \to 0,$ 
\begin{equation*}
\dtv (\nu^\epsilon, \mu^\epsilon) \to 1. 
\end{equation*} 
\end{lem}  
Note that this lemma and Pinsker's inequality \eqref{pinskerine} imply that the Kullback-Leibler divergence remains at least order one in the small $\epsilon$ limit, unless $m$ solves \eqref{ODE}.

\section{Gaussian Approximation of SDEs via ODEs}\label{sec:cts-cts}
Let $\nu^\epsilon$ be the law of $l^\epsilon$ defined via the linear SDE
\begin{equation}\label{approxeqcont}
dl^\epsilon(t)= \Bigl( f\bigl(m(t)\bigr) + Df\bigl(m(t)\bigr)\bigl(l^\epsilon(t)-m(t)\bigr)\Bigr)dt + \sqrt{\epsilon\Sigma} \, dW, \quad l^\epsilon(0)\sim N\bigl(m_0,\epsilon C_0\bigr).
\end{equation}
The main result of this section is  Theorem  \ref{theorem:cts} below. It shows that, for  small $\epsilon$, the Gaussian measure $\nu^\epsilon$  accurately approximates  $\mu^\epsilon$, the law of $v^\epsilon$ given by \eqref{sde}.
The proof is based on two observations. First, that $\nu^\epsilon$ has time-marginals $\nu_t^\epsilon = N\bigl(m(t),C^\epsilon(t)\bigr)$, where $m$ solves the ODE \eqref{ODE}, and $C^\epsilon=\epsilon C$ with $C$ solving \eqref{ODEcov}. 
Second, the following lemma:

\begin{lem}\label{lemmaconttime}
Let $\nu_0^\epsilon \ll \mu_0^\epsilon,$ let $g_t: \R^D\to \R^D,$ and let $\mu^\epsilon$ and $\nu^\epsilon$ be the laws in $C([0,T],\R)$ of $v^\epsilon$ and $\l^\epsilon$ given by
\begin{align*}
dv^\epsilon &= f(v^\epsilon)dt +  \sqrt{\epsilon\Sigma} \, dW, \quad v(0)\sim \mu_0^\epsilon, \\
dl^\epsilon &= g_t(l^\epsilon)dt +  \sqrt{\epsilon\Sigma} \, dW, \quad l^\epsilon(0)\sim \nu_0^\epsilon.
\end{align*}
Then
\begin{equation*}
\dkl(\nu^\epsilon\| \mu^\epsilon) =\dkl(\nu_0^\epsilon\|\mu_0^\epsilon)+ \frac{1}{2\epsilon} \E^{\nu^\epsilon} \int_0^T |f(v)-g_t(v)|_\Sigma^2\,dt.
\end{equation*}
\end{lem}
\begin{proof}
Girsanov's theorem \cite{elworthy} gives
\begin{equation*}
\frac{d\nu^\epsilon}{d\mu^\epsilon}(v)= \frac{d\nu_0^\epsilon}{d\mu_0^\epsilon}\bigl(v(0)\bigr) \exp\Bigl(\frac{1}{\epsilon} \int_0^T\langle g_t(v)-f(v), dv\rangle_\Sigma -\frac{1}{2\epsilon}\int_0^T |g_t(v)|_{\Sigma}^2-|f(v)|^2_\Sigma\, dt\Bigr).
\end{equation*}
Therefore,
\begin{align*}
\dkl(\nu^\epsilon\|\mu^\epsilon)&= \E^{\nu^\epsilon} \log\Bigl(\frac{d\nu^\epsilon}{d\mu^\epsilon}\Bigr)\\
&=\dkl(\nu_0^\epsilon\|\mu_0^\epsilon) + I,
\end{align*}
where using the martingale property of Ito's integral
\begin{align*}
I&=\frac{1}{\epsilon} \E^{\nu^\epsilon} \Bigl(\int_0^T\langle g_t(v)-f(v), dv\rangle_\Sigma -\frac{1}{2}\int_0^T |g_t(v)|_{\Sigma}^2-|f(v)|^2_\Sigma\, dt \Bigr)\\
&=\frac{1}{\epsilon}\E^{\nu^\epsilon}\Bigl( \int_0^T \langle g_t(v)-f(v) ,\sqrt{\epsilon \Sigma}\, dW \rangle_\Sigma+\frac{1}{2} \int_0^T |f(v)-g_t(v)|_\Sigma^2\,dt\Bigr)\\
&= \frac{1}{2\epsilon} \E^{\nu^\epsilon} \int_0^T |f(v)-g_t(v)|_\Sigma^2\,dt.
\end{align*}
\end{proof}

\begin{thm}\label{theorem:cts}
Suppose that Assumption \ref{assumptionode} holds. Let $\mu^\epsilon$ and $\nu^\epsilon$ be the laws in $C([0,T],\R^D)$ of, respectively, $v^\epsilon$ given by \eqref{sde} and  $l^\epsilon$ given by \eqref{approxeqcont}. 
Then there is $c>0$ such that, for all $\epsilon$ sufficiently small,
$$\dkl(\nu^\epsilon\| \mu^\epsilon)\le \dkl(\nu_0^\epsilon\|\mu_0^\epsilon) +  c \epsilon.$$
Moreover, for $t\in [0,T],$ 
$$\dkl(\nu_t^\epsilon\| \mu_t^\epsilon)\le \dkl(\nu_0^\epsilon\|\mu_0^\epsilon) + c\epsilon,$$
where $\nu_t^\epsilon$ and $\mu_t^\epsilon$ denote the time-marginals of $\nu^\epsilon$ and $\mu^\epsilon$ at time $t.$
\end{thm}
\begin{proof}
If $\nu_0^{\epsilon}$ is not absolutely continuous with respect to $\mu_0^{\epsilon},$ then $\dkl(\nu_0^\epsilon\|\mu_0^\epsilon)=\infty$  and the result is trivial. Thus we assume that absolute continuity holds.
Throughout the proof $c$ denotes a positive constant which is independent of all sufficiently small $\epsilon,$ and may vary from line to line. As noted before, $\nu^\epsilon$ has time-marginals $\nu_t = N\bigl(m(t),C^\epsilon(t)\bigr)$ with $C^\epsilon=\epsilon C$. This, combined with Lemma \ref{lemmaconttime} with the choice $g_t(\cdot):= f\bigl(m(t)\bigr) + Df\bigl(m(t)\bigr)\bigl(\cdot-m(t)\bigr)$, gives
\begin{align}\label{intermediatestep}
\begin{split}
\dkl(\nu^\epsilon\| \mu^\epsilon)&=  \dkl(\nu_0^\epsilon\|\mu_0^\epsilon)  + \frac{1}{2\epsilon}\int_0^T \E^{N(m(t),C^\epsilon(t))}|f(v)-g_t(v)|_{\Sigma}^2 dt.
\end{split}
\end{align}
Recall Taylor's formula with reminder
\begin{align*}
\bigl(f(u)-g_t(u)\bigr)_d= 2 \sum_{|\alpha|=2}\frac{(u-m)^\alpha}{\alpha!}\int_0^1 (1-t)^2 \partial^\alpha f_d\bigl(m+t(u-m)\bigr)\, dt, \quad 1\le d\le D,
\end{align*}
to deduce, using Assumption \ref{assumptionode}, that
\begin{align*}
|\bigl(f(u)-g_t(u)\bigr)|&\le c|u-m|^2 \int_0^1 1+|m+t(u-m)|^s \,dt \\
&\le c\sum_{r=0}^s |u-m|^{r+2}. 
\end{align*}
Thus, 
\begin{equation}\label{taylor}
|f(u)-g_t(u)|^2 \le c \sum_{r=0}^s |u-m|^{2r+4}.
\end{equation}
Combining \eqref{intermediatestep} with \eqref{taylor} and using that $C^\epsilon = \epsilon C$ yields, for all $\epsilon$ sufficiently small,
\begin{align*}
\begin{split}
\dkl(\nu^\epsilon\| \mu^\epsilon)&=  \dkl(\nu_0^\epsilon\|\mu_0^\epsilon)  + \frac{1}{2\epsilon}\int_0^T \E^{N(m(t),C^\epsilon(t))}|f(v)-g_t(v)|_{\Sigma}^2 dt \\
&\le \dkl(\nu_0^\epsilon\|\mu_0^\epsilon)  +  \frac{c}{\epsilon}\int_0^T \E^{N(m(t),C^\epsilon(t))}\sum_{r=0}^s |v(t)-m(t)|^{2r+4} dt \\
&\le \dkl(\nu_0^\epsilon\|\mu_0^\epsilon)  +  \frac{c}{\epsilon}\max_{0\le r \le s} \int_0^T \Bigl( \E^{N(m(t),C^\epsilon(t))} |v(t)-m(t)|^{2}\Bigr)^{r+2} dt.\\
&\le \dkl(\nu_0^\epsilon\|\mu_0^\epsilon)  +  c\epsilon,
\end{split}
\end{align*}
 which completes the proof of the first claim. The bound for the marginals then follows from \eqref{marginalsbound}.
\end{proof}

\begin{remark}
It is not difficult to see that if $f$ is linear then 
\begin{equation*}
\dkl(\nu^\epsilon\| \mu^\epsilon) = \dkl(\nu_0^\epsilon\| \mu_0^\epsilon).
\end{equation*}
This is the well-known Kullback-Leibler stability of the Fokker-Planck equation with respect to initial conditions, of which a  more general version can be found in \cite{clark1999relative}.
\end{remark}

\section{Gaussian Approximations of SDEs via discretized ODEs}\label{sec:cts-disc}\
In this section we study the approximation of $\mu^\epsilon$, the law in $C([0,T],\R^D)$ of $v^\epsilon$ given by \eqref{sde}, by a Gaussian $\nu^{\epsilon, \Delta t}$. The path measure $\nu^{\epsilon, \Delta t}$ will be constructed so that the means and covariances of the time-marginals are given by numerical approximations with step size $\Delta t$ to the ODEs \eqref{ODE} and \eqref{ODEcov}, respectively. For simplicity, we study the effect of discretizing the ODEs with 
the Euler method, but the results extend with no effort to other Runge-Kutta methods and numerical schemes. Theorem \ref{theorem:ctsdisc} below bounds $\dkl(\nu^{\epsilon, \Delta t}\|\mu^\epsilon)$ in terms of $\epsilon$ and the step-size $\Delta t$.

We now spell out the construction of the measure $\nu^{\epsilon,\Delta t}.$
Given an integer $K>0,$ let $\Delta t = T/K$ and define, for $1\le k \le K,$
\begin{align}\label{eq:discretizationsec3}
\begin{split}
m_{k+1}^{\Delta t} &:= m_k^{\Delta t}  +\Delta t  f(m_k^{\Delta t}), \quad m_0^{\Delta t}=m_0,\\
C_{k+1}^{\Delta t} &:=C_k^{\Delta t} + \Delta t\Bigl( Df(m_k^{\Delta t})C_k^{\Delta t}  +  C_k^{\Delta t} Df(m_k^{\Delta t})^T +  \Sigma\Bigr), \quad C_0^{\Delta t} =C_0.
\end{split}
\end{align}

Let $t_k=k\Delta t,\, 0\le k \le K, $  and define piecewise linear functions $m^{\Delta t}$ and $C^{\Delta t}$ in $[0,T]$ by interpolation. That is, for $t\in (t_k,t_{k+1})$,

\begin{align}\label{eq:timemarginals}
\begin{split}
m^{\Delta t}(t) &:= m_k^{\Delta t} + (t-t_k)f(m_k^{\Delta t}),\\
C^{\Delta t}(t) &:= C_k^{\Delta t} + (t-t_k)\Bigl( Df(m_k^{\Delta t})C_k^{\Delta t}  +  C_k^{\Delta t} Df(m_k^{\Delta t})^T + \Sigma\Bigr).
\end{split}
\end{align}

Finally, we let $\nu^{\epsilon,\Delta t}$ be the law of $l^{\epsilon,\Delta t}$ defined via the piecewise linear SDE
\begin{equation}
\label{eq:plsde}
dl^{\epsilon, \Delta t} = g_t^{\Delta t} (l^{\epsilon, \Delta t}) \, dt + \sqrt{\epsilon\Sigma} \, dW, \quad l^{\epsilon, \Delta t}(0)\sim N\bigl(m_0,\epsilon C_0\bigr),
\end{equation}
with
\begin{equation}\label{eq:driftapproxode}
g_t^{\Delta t}(\cdot):= f(m_k^{\Delta t}) + Df(m_k^{\Delta t})\bigl( \,\,\cdot \, - m^{\Delta t}(t)\bigr),\quad  t\in(t_k,t_{k+1}).
\end{equation}
By construction $\nu^{\epsilon, \Delta t}$ has time marginals $\nu_t^{\epsilon, \Delta t}=N\bigl(m^{\Delta t}(t), C^{\epsilon,\Delta t}(t)\bigr)$ with $m^{\Delta t}$ and $C^{\epsilon, \Delta t}:=\epsilon C^{ \Delta t}$ defined by \eqref{eq:timemarginals}.

\begin{remark}
For fixed $\Delta t>0$, $m^{\Delta t}\neq m$ except in trivial cases. Thus, by Lemma \ref{lemmaldp}  it is necessary to let $\Delta t$ depend on $\epsilon$  in order to have accurate approximations $\nu^{\epsilon, \Delta t}$ of $\mu^\epsilon$ in the limit $\epsilon \to 0$.
\end{remark}

\begin{thm}\label{theorem:ctsdisc}
Suppose that Assumption \ref{assumptionode} holds. Let $\mu^\epsilon$ be the law in $C([0,T],\R^D)$ of $v^\epsilon$ given by \eqref{sde}, and let $\nu^{\epsilon, \Delta t}$ be as above. 
For $K\in \N$ let $\Delta t = T/K.$ Then there is $c>0,$ independent of all sufficiently small $\Delta t$ and $\epsilon,$ such that
$$\dkl(\nu^{\epsilon,\Delta t}\| \mu^\epsilon)\le \dkl(\nu_0^{\epsilon,\Delta t}\|\mu_0^\epsilon) +  c \epsilon +c \frac{(\Delta t)^2}{\epsilon}.$$
Moreover, for $t\in (0,T],$ 
$$\dkl(\nu_t^{\epsilon,\Delta t}\| \mu_t^\epsilon)\le \dkl(\nu_0^{\epsilon,\Delta t}\|\mu_0^\epsilon) +  c \epsilon +c \frac{(\Delta t)^2}{\epsilon},$$
where $\nu_t^{\epsilon,\Delta t}$ and $\mu_t^{\epsilon,\Delta t}$ denote the marginals of $\nu^{\epsilon, \Delta t}$ and $\mu^{\epsilon, \Delta t}$ at time $t.$
\end{thm}
\begin{proof}
Throughout the proof $c$ is a positive constant that may change from line to line, and is independent of all sufficiently small $\epsilon$ and $\Delta t$. We recall that  $\nu^{\epsilon,\Delta t}=N(m^{\Delta t},C^{\epsilon,\Delta t})$ is the law of the SDE \eqref{eq:plsde}. This, combined with Lemma \ref{lemmaconttime} applied with $g_t = g_t^{\Delta t}$ as defined by \eqref{eq:driftapproxode}, gives
\begin{align}\label{intermediatestep2}
\begin{split}
\dkl&(\nu^{\epsilon, \Delta t}\| \mu^\epsilon)=  \dkl(\nu_0^{\epsilon,\Delta t}\|\mu_0^\epsilon)  + \frac{1}{2\epsilon}\int_0^T \E^{\nu_t^{\epsilon, \Delta t}}|f(v)-g_t^{\Delta t}(v)|_{\Sigma}^2\,\, dt\, .
\end{split}
\end{align}
We split the integral as follows:
\begin{equation}
|f(v)-g_t^{\Delta t}(v)|_{\Sigma}^2 \le 
c|g_t(v) - g_t^{\Delta t}(v)|^2_{\Sigma}
+c|f(v) - g_t(v)|^2_\Sigma, 
\end{equation}
and bound each of the two terms.
For the first one note that, for $t \in (t_k, t_{k+1})$ and all sufficiently small $\Delta t$,
\begin{align*}
|g_t(v) - g_t^{\Delta t}(v)| &\le |f\bigl (m(t)\bigr) - f(m_k) | + |Df(m_k)| |m(t)-m_k| +\cdots \\& |Df\bigl(m(t)\bigr)- Df(m_k)| |v(t)-m^{\Delta t}(t)| + 
 |Df\bigl(m(t)\bigr)- Df(m_k)| |m^{\Delta t}-m_k|\\ 
& \le \,c\Delta t + c\Delta t  |v(t)-m^{\Delta t}|.
\end{align*}
Thus, 
\begin{equation}
 |g_t(v) - g_t^{\Delta t}(v)|^2_\Sigma \le \,c (\Delta t)^2 + c(\Delta t)^2  |v(t)-m^{\Delta t}|^2.
\end{equation}
For the second one, as in the proof of Theorem \ref{theorem:cts}, we have for small enough $\Delta t$ that
\begin{align*}
|f(v)- g_t(v)|^2 &\le \,c\sum_{r=0}^s |v-m|^{r+2}\\
&\le \,c\sum_{r=0}^s |v-m^{\Delta t}|^{r+2} + c\sum_{r=0}^s |m^{\Delta t} -m|^{r+2}\\
&\le \,c\sum_{r=0}^s |v-m^{\Delta t}|^{r+2}  + c(\Delta t)^2.
\end{align*}
Putting everything together, for all $\Delta t$ and $\epsilon$ sufficiently small,
\begin{align*}
\begin{split}
\dkl&(\nu^{\epsilon, \Delta t}\| \mu^\epsilon)=  \dkl(\nu_0^{\epsilon,\Delta t}\|\mu_0^\epsilon)  + \frac{1}{2\epsilon}\int_0^T \E^{\nu_t^{\epsilon, \Delta t}}|f(v)-g_t^{\Delta t}(v)|_{\Sigma}^2 \,\,dt \\
&\le \dkl(\nu_0^{\epsilon,\Delta t}\|\mu_0^\epsilon)  +  \frac{c}{\epsilon} (\Delta t)^2 + \frac{c}{\epsilon} \int_0^T  \E^{N(m^{\Delta t}(t),C^{\epsilon, \Delta t})}  \sum_{r=0}^s |v-m^{\Delta t}|^{r+2} \\
&\le \dkl(\nu_0^{\epsilon,\Delta t}\|\mu_0^\epsilon) +  \frac{c}{\epsilon} (\Delta t)^2 + c \epsilon.
\end{split}
\end{align*}
\end{proof}

\section{Gaussian Approximation of Discretized SDEs via Discretized ODEs}\label{sec:disc-disc}

Consider the Euler-Maruyama discretization of \eqref{sde}
\begin{equation}\label{sdeeuler}
v_{k+1}^\epsilon= v_k^\epsilon + f(v_k^\epsilon) \Delta t +  \sqrt{\epsilon \Sigma \Delta t}\, \xi_k, \quad v_0^\epsilon\sim \mu_0^\epsilon,
\end{equation}
where the $\xi_k$ are independently drawn from a standard Gaussian distribution. In this section we consider $\Delta t>0$  small and fixed, and analyze small $\epsilon$ limits. A structure of the form \eqref{sdeeuler} arises from discretization of SDEs, but also in many stochastic algorithms \cite{kushner}. In the case of interest where the function $f$ is nonlinear, the distribution $\mu_k^\epsilon$ of $v_k^\epsilon$ is not Gaussian. We study the approximation of these measures by Gaussians $\nu_k^\epsilon$, whose means and covariances are built using discretizations of the ODEs \eqref{ODE} and \eqref{ODEcov}. 

We now detail the construction of the measures $\nu_k^\epsilon$. 
Let 
\begin{align}
m_{k+1} &= m_k + f(m_k)\Delta t, \label{eq:eulermeanODE} \\
C_{k+1} &= \bigl(I+Df(m_k)\Delta t\bigr) C_k \bigl(I+Df(m_k)\Delta t\bigr)^T+ \Sigma\Delta t. \label{eq:disccovODE}
\end{align}
These agree with the discretization \eqref{eq:discretizationsec3} used in Section \ref{sec:cts-disc}, except for an extra  $(\Delta t)^2$ term in the covariance.
We set  $C_k^\epsilon := \epsilon C_k$, and finally $\nu_k^\epsilon=N(m_k, C_k^\epsilon).$

The subsequent analysis is parallel to that of the previous sections. We again use two observations. First, that the $\nu_k^\epsilon$ are the laws of $l_k^\epsilon$ given by 
\begin{equation}\label{linearcomplete}
l_{k+1}^\epsilon = l_k^\epsilon + \bigl( f(m_k) + Df(m_k)(l_k^\epsilon-m_k)  \bigr) \Delta t + \sqrt{\epsilon\Sigma\Delta t}\,\, \xi_k.
\end{equation}
Second, the following lemma, analogous to Lemma \ref{lemmaconttime}. A derivation of this result can be found for instance in the appendix of \cite{archambeau}. We include a short proof that highlights how the chain rule Lemma \ref{lemmakl} makes the Kullback-Leibler divergence  well suited for the analysis of conditionally Gaussian dynamics, as those defined by \eqref{sdeeuler}.
\begin{lem}\label{lemmageneralbounddiscretetime}
Let $g_k:\R^D\to \R^D$. Let $\mu_k^\epsilon$ and $\nu_k^\epsilon$ be the laws of $v_k^\epsilon$ and $l_k^\epsilon$ in $\R^D$ given by
\begin{align*}\label{sdeeulerlinear}
l_{k+1}^\epsilon&= l_k^\epsilon + g_k(l_k^\epsilon) \Delta t +  \sqrt{\epsilon\Sigma\Delta t} \,\, \xi_k, \quad l_0\sim \nu_0^\epsilon, \\
v_{k+1}^\epsilon&= v_k^\epsilon + f(v_k^\epsilon) \Delta t +  \sqrt{\epsilon \Sigma \Delta t}\,\, \xi_k, \quad v_0^\epsilon\sim \mu_0^\epsilon.
\end{align*}
Denote by $ \mu_{0:k}^\epsilon$ and $ \nu_{0:k}^\epsilon$ the law of $(v_j^\epsilon)_{j=0}^k$ and $(l_j^\epsilon)_{j=0}^k$  in $\R^{D (k+1)}$. Then 
\begin{equation}\label{eq:formulakldiscrete}
\dkl (\nu_{0:k}^\epsilon \|\mu_{0:k}^\epsilon) = \dkl( \nu_0^\epsilon \| \mu_0^\epsilon) + \frac{\Delta t}{2\epsilon} \sum_{j=0}^{k-1} \E^{\nu_j^\epsilon} |f-g_j|_\Sigma^2.
\end{equation}
\end{lem}
\begin{proof}
We show that 
\begin{equation}\label{eq:inductionstep}
\dkl (\nu_{0:k+1}^\epsilon \|\mu_{0:k+1}^\epsilon) =\dkl (\nu_{0:k}^\epsilon \|\mu_{0:k}^\epsilon) + \frac{\Delta t}{2\epsilon}\E^{\nu_k^\epsilon} |f-g_k|_\Sigma^2.
\end{equation}
Iterating \eqref{eq:inductionstep} gives \eqref{eq:formulakldiscrete}.
To prove \eqref{eq:inductionstep} we apply Lemma \ref{lemmakl} with $\nu=(\nu_{0:k},\nu_{k+1}),\, \mu=(\mu_{0:k},\mu_{k+1})$, and the kernels
\begin{equation*}
\mu(\cdot|x) = N\Bigl(x + f(x)\Delta t, \epsilon \Sigma \Delta t\Bigr), \quad \nu(\cdot|x)= N\Bigl(x + g_k(x)\Delta t, \epsilon \Sigma \Delta t\Bigr).
\end{equation*}
To conclude, recall that the Kullback Leibler between two Gaussians with the same covariance is
\begin{equation*}\label{KLgaussians}
\dkl\bigl(N(m,C)\|N(\bar{m}, C)\bigr) =\frac12  |\bar{m} - m|^2_{C}.
\end{equation*}
\end{proof}

The following result shows that the Gaussians $\nu_k^\epsilon$ accurately approximate the distribution $\mu_k^\epsilon$ of $v_k^\epsilon$ in the small noise limit. More so, the result holds in path space.

\begin{thm}\label{theoremdiscrete}
Suppose that Assumption \ref{assumptionode} holds. Let $\mu_{0:k}^\epsilon$ be the law of $(v_j^\epsilon)_{j=0}^k$ given by \eqref{sdeeuler}, and let $\nu_{0:k}^\epsilon$ be the law of $(l_j^\epsilon)_{j=0}^k$ given by \eqref{linearcomplete}. Then, for all sufficienly small $\Delta t$ and $\epsilon$,  and for all $k$ with $k\Delta t<T,$ 
\begin{equation*}
\dkl (\nu_{0:k}^\epsilon \|\mu_{0:k}^\epsilon) \le \dkl(\nu_0^\epsilon\| \mu_0^\epsilon) + c\epsilon, \quad \dkl (\nu_k^\epsilon \|\mu_k^\epsilon) \le \dkl(\nu_0^\epsilon\| \mu_0^\epsilon) + c\epsilon,
\end{equation*}
where $c>0$ is independent of the noise strength $\epsilon.$ 
\end{thm}
\begin{proof}
Throughout the proof $c$ denotes a positive constant that is independent of $\epsilon$ and may change from line to line.
We use the previous lemma with the choice
\begin{equation*}
g_k(u):=f(m_k) + Df(m_k)(u-m_k)
\end{equation*}
together with the observation made in \eqref{linearcomplete}
to deduce that, for $k$ with $k\Delta t\le T,$ 
\begin{align}\label{eq:partial1}
\begin{split}
\dkl (\nu_{0:k}^\epsilon \|\mu_{0:k}^\epsilon) &=  \dkl(\nu_0^\epsilon\| \mu_0^\epsilon) + \frac{\Delta t}{2\epsilon}\sum_{j=0}^{k-1}  \E^{\nu_j^\epsilon} |f-g_j|_\Sigma^2\\
 &\le  \dkl(\nu_0^\epsilon\| \mu_0^\epsilon) +\frac{c}{\epsilon}\max_{k:k\Delta t\le T} \E^{N(m_k,C_k^\epsilon)}|f(u)-g_k(u)|^2.
 \end{split}
\end{align}

Now, using Remark \ref{rem:disc-cov-bound} it can be shown as in Theorem \ref{theorem:cts} that 
\begin{equation}\label{eq:partial2}
|f(u)-g_k(u)|^2 \le c \sum_{r=0}^s |u-m_k|^{2r+4}.
\end{equation}

To conclude we combine  \eqref{eq:partial1} and \eqref{eq:partial2},  and recall that $C_k^\epsilon=\epsilon C_k$ and Remark \ref{rem:disc-cov-bound} to deduce that, for $k$ with $k\Delta t\le T$ and all sufficiently small $\epsilon,$  
\begin{align*}
\dkl (\nu_{0:k}^\epsilon \|\mu_{0:k}^\epsilon)
&\le  \dkl(\nu_0^\epsilon\| \mu_0^\epsilon) + \frac{c}{\epsilon}\max_{k:k\Delta t\le T}  \E^{N(m_k,C_k^\epsilon)}\sum_{r=0}^s |u-m_k|^{2r+4}\\
&\le \dkl(\nu_0^\epsilon\| \mu_0^\epsilon)  +c\epsilon.
\end{align*}
\end{proof}

\begin{example}
Again if $f$ is linear then, for $k\ge 1,$
\begin{equation*}
\dkl(\nu_{0:k}^\epsilon\| \mu_{0:k}^{\epsilon})= \dkl(\nu_0^\epsilon\| \mu_0^\epsilon), \quad \dkl(\nu_k^\epsilon\| \mu_k^{\epsilon})\le \dkl(\nu_0^\epsilon\| \mu_0^\epsilon).
\end{equation*} 
As an example, suppose that $D=1,$  $f(u)=au+b,$ and $\mu_0^\epsilon$ is Gaussian in $\R.$ Then a direct proof of the last inequality can be easily obtained. Indeed, it boils down to showing that, for any $m,\tilde{m}, \sigma,\tilde{\sigma}, a,b\in \R$ and $\epsilon>0,$
\begin{equation*}
\dkl\bigl(N(am+b,a^2\sigma^2+\epsilon)\| N(a\tilde{m}+b,a^2\tilde{\sigma}^2+\epsilon)\bigr)\le \dkl\bigl(N(m,\sigma^2)\| N(\tilde{m},\tilde{\sigma}^2)\bigr).
\end{equation*} 
In other words, 
\begin{equation*}
\frac{a^2\sigma^2 + \epsilon}{a^2\tilde{\sigma}^2 + \epsilon} -\log\Bigl(\frac{a^2\sigma^2+\epsilon}{a^2\tilde{\sigma}^2+\epsilon}\Bigr) -1 +\frac{a^2(m-\tilde{m})^2}{a^2\tilde{\sigma}^2+\epsilon} \le \frac{\sigma^2}{\tilde{\sigma}^2}-\log\Bigl(\frac{\sigma^2}{\tilde{\sigma}^2}\Bigr) -1 + \frac{(m-\tilde{m})^2}{\tilde{\sigma}^2},
\end{equation*}
which is immediate since the sum of the first two terms in the LHS is smaller than the sum of the first two terms in the RHS and the last term in the LHS is smaller than the last term in the RHS.
\end{example}

          %







          %


          \end{document}